\documentclass[a4paper,10pt]{article}
\usepackage[toc,page]{appendix}
\usepackage{amsfonts}
\usepackage{mathrsfs}
\usepackage{amsmath}
\usepackage{amssymb}
\usepackage{verbatim}
\usepackage{epsfig}
\usepackage{cite}
\usepackage{amstext}
\usepackage{tabularx,booktabs}
\usepackage{graphicx}
\usepackage{subfig}
\usepackage{color}
\usepackage{hyperref}
\usepackage{esint}
\usepackage{authblk}
\usepackage[toc,page]{appendix}
\usepackage{blkarray}
\usepackage{multirow}
\hypersetup{
    colorlinks,
    citecolor=black,
    filecolor=black,
    linkcolor=black,
    urlcolor=black
}
\numberwithin{equation}{section}
\numberwithin{figure}{section}

\newcolumntype{C}{>{\centering\arraybackslash}X} 
\usepackage{lipsum}

\begin{document}
\title{A 3D Nonlinear Maxwell's Equations Solver Based On A Hybrid Numerical Method}
\author{Aihua Lin}
\author{Per Kristen Jakobsen}
\affil{Department of Mathematics and Statistics, UIT the Arctic University of Norway, 9019 Troms\o, Norway}

\renewcommand\Authands{ and }
\maketitle

\begin{abstract}
In this paper we explore the possibility for solving the 3D Maxwell's equations in the presence of nonlinear and/or inhomogeneous material response.  We propose using a hybrid approach which combines a boundary integral representation with a domain-based method. This hybrid approach has previously been successfully applied to 1D linear and nonlinear transient wave scattering problems. The basic idea of the approach is to propagate the Maxwell's equations inside the scattering objects forward in time by using a domain-based  method, while a boundary integral representation of the electromagnetic field is used to supply the domain-based method with the required surface values. Thus no grids outside the scattering objects are needed and this greatly reduces the computational cost and complexity. 

\end{abstract}

\section{\bigskip Introduction}
Boundary Element method (BEM), as a tool for solving scattering problems, has several attractive features.
First and foremost, BEM is well suited to treating scattering problems in unbounded domains because the boundary integral equations are located on the surfaces of the scattering objects and thus one whole dimension is taken out of the problem. Also, the scattering objects are usually defined by sharp material boundaries and thus a domain-based method must seek to resolve the fast variation in the corresponding solutions generated by the boundaries. This problem has been more or less solved by using perfectly matched layers\cite{PML1,PML2}, but the solution comes with additional cost and complexity.

 Secondly,  the surface localization of boundary integral equations also means that the boundary discretization, which leads to the BEM equations, can be optimized with respect to the geometry of each surface separately when there are more than one scattering object, which there usually is. For domain-based methods like Finite Element method (FEM)\cite{FEM1,FEM2,FEM3} and Finite Difference Time Domain method (FDTD), such an optimization can only be achieved by using non-uniform or multiple grids tailored to the geometrical shape of the objects. These kinds of efforts have met with some degree of success, but it does add new layers of complexity. 
 
 Thirdly,  BEM,  is exceptionally well suited for modeling scattering problems where the sources are slow on relevant timescales. In this setting the boundary integral equations are derived in the spectral domain and the discretized equations defining BEM  only needs to be solved for the small set of discrete frequencies that are required for an accurate representation of the time dependent source. For a domain-based method like FDTD, near-stationary sources are the worst possible case since FDTD is, as the name indicate, a time domain method and slow sources mean that the Maxwell's equations  have to be solved for a long interval of time, which is costly in terms of computational resources. 
 
 Given all these attractive features of BEM for solving scattering problems in electromagnetics, it is somewhat surprising that in a popularity contest, FEM and FDTD beat BEM hands down\cite{Cheng}. The reason for this is that in addition to all its attractive features, BEM has some real drawbacks too.
 
 Firstly, reformulating the scattering problem for a system of PDEs in terms of boundary integral equations, typically will involve classical but fairly intricate mathematical tools, as will the discretization of the resulting integral equations. In particular, one will need to content with singularities which appear in the limits that are always taken while deriving the boundary equations from the PDEs. Accurately representing these singularities in the ensuing discretization leading to BEM, is a major issue. FEM has some of the same issues, whereas FDTD is much simpler to implement using mathematical tools that are straight forward and known to all.
 
Secondly, BEM relies on Green's functions,  and such functions can only be defined for linear systems of PDEs. In our opinion this is the major reason why BEM is significantly less popular than the major contenders FEM and FDTD. Nonlinearities are common in most areas of pure and applied science, and having to change your whole computational approach when nonlinearities are added to your model is a major nuisance. In some areas of application nonlinearities dominate and hence a computational approach based on Green's functions is out of the question, fluid mechanics is such an area. However in other areas one can frequently disregard nonlinearities and a computational approach based on Green's functions is feasible. Scattering of electromagnetic waves, where nonlinearities only come into play at very high field intensities, is such an area.

Thirdly, BEM, which is simplest to formulate in the spectral domain, is at it's best when we are looking at a near stationary  situation with narrow band sources. For transient scattering with broad band sources, one needs to formulate BEM in the time domain, and this involves space-time Green's functions that in general tends to be more singular than the frequency domain ones. However, for the specific case of electromagnetic scattering an integral formulation has been derived \cite{Jones64} whose singularities are fairly weak. This integral formulation of electromagnetics is the foundation of all time dependent BEM formulations to date, and is also the basis for our approach. However, the progress in developing these time dependent BEM schemes has been slow due to several drawbacks.

  The time-domain integral equations are retarded and this means that in order to compute 
the solutions at a certain time, one needs to retain the solutions for an interval of previous times that in some cases can be very large. This leads to a large memory requirement that needs to be met using parallel processing. In today's computational environment parallel  implementation of  time dependent BEM is fairly standard, but the possibly limited efficiency due to the problem of load balancing is something that always must be contended with.

However the major obstacle that has prevented this method from being widely applied for electromagnetic scattering  is the occurrence of numerical instabilities. These instabilities, whose source is not fully understood, occur not at early times, but at later times, and have become known as the late time instability. Many efforts have been made and several techniques have been developed in order to  improve the stabilities of BEM schemes for time dependent electromagnetics  in the last several decades \cite{Weile, Weile2, Shanker, Walker, Walker2,Zhao, Huang}. 
Our work has not been aimed at joining in or improving on any of the efforts pursued by these research groups. Our major aim has been to generalize boundary approaches to electromagnetic scattering for cases where the scattering objects has an inhomogeneous and/or nonlinear response. Most work in the area of time dependent BEM has been focused on metallic objects with linear response whose dimensions are large with respect to the wave length, antennas is a major example of the kind of structure one has been interested in. The kind of problem we have in mind is scattering from very small objects, from micron to nanometer scale, objects who might have engineered inhomogeneities  in their structure and strong, also engineered, nonlinear response. At this scale, the standard simplifying assumption of disregarding the inside of the scattering objects and modeling them using  surface charges and currents, is not applicable. The skin depth at this microscopic scale can easily be as large as the scattering objects themselves, which is a marked difference to what is true for macroscopic antenna theory. For this reason, and also because of the complex inner structure of these microscopic scattering objects, a different approach is needed. The traditional BEM can not be used here. 

Our approach is based on the same integral representation of the electromagnetic field \cite{Jones64} as the traditional BEM, however, we use the integral representation in a way that is different  from what one does in BEM.  We solve the Maxwell's equations on the inside of each scattering object, as an initial boundary value problem, and use the integral identities to supply the boundary values needed in order to make the initial boundary value problem for Maxwell, well posed. This kind of approach for solving electromagnetic scattering problems was first proposed in 1972 by E. Wolf and D. N. Pattanyak \cite{wolf} in the context of stationary linear scattering and was based on the Ewald-Oseen optical extinction theorem. For this reason we call this particular way of reformulating the electromagnetic scattering problem for the Ewald Oseen Scattering(EOS) formulation. This reformulation can be applied to any kind of wave scattering situation. It has previously been applied to two toy models of 1D linear and nonlinear transient wave scattering by the authors \cite{Aihua} where the EOS formulations work perfectly well with high accuracy and low computational load and without any instabilities, even at very late times.

In section  \ref{EOS3d} we show some of the details of the derivation of our EOS formulation for Maxwell's equations and in section \ref{implementation3d} we discuss some tests we have run on our numerical implementation of the EOS formulation of Maxwell's equations.
In this paper we do not describe the numerical implementational details, like how we handle the singular integrals and issues of numerical stability. These, mainly very technical considerations, would cloud the main message of the current paper, which is that our EOS formulation of scattering problems works. The technical details pertaining to our choise of implementation, some of which are probably relevant for most numerical implementations of the EOS formulation, will be reported elsewhere in a paper soon to appear.  Here we will just note that, just like for the 1D case, the internal numerical scheme, Lax-Wendroff for our case, determines a stability interval for the time step. The difference is that,  in the 1D case, the stability interval  is purely determined by the internal numerical scheme while in 3D case, there is another lower limit of the stability interval determined by the integral part of the scheme.  We also find that  the late time instability is highly depended on the features of the scattering materials.  Section \ref{conclusion3d} summarizes what we have achieved and discuss extensions of our work that could be of interest to pursue. 
\section{EOS formulations of the 3D Maxwell's equations}\label{EOS3d}
In this paper, we investigate an electromagnetic scattering problem described by the 3D Maxwell's equations 
\begin{equation*}
\begin{split}
\nabla \times {\bf E}+\partial_t {\bf B}&=0,\\
\nabla \times {\bf H}-\partial_t {\bf D}&={\bf J},\\
\nabla \cdot {\bf D}=\rho,\\
\nabla \cdot {\bf B}=0,
\end{split}
\end{equation*}
where ${\bf J}$ and $\rho$ are the current density and the charge density of free charges. Bound charges and currents determine ${\bf D}$ and ${\bf H}$ as functionals of ${\bf E}$ and ${\bf B}$,
\begin{equation*}
\begin{split}
{\bf D}&={\bf D}[\bf E,\bf B],\\
\bf H&=\bf H[\bf E,\bf B].
\end{split}
\end{equation*}
In the simplest situation, where the response from the bound charges
and currents is linear, isotropic, homogeneous and instantaneous,  we
have
$$
{\bf D}=\varepsilon {\bf E},\ {\bf H}=\frac{1}{\mu}\bf B,
$$
where $\sqrt{\frac{1}{\varepsilon \mu}}=c$ is the speed of light in the material. For this particular situation, we have
\begin{equation}\label{eq:maxwell}
\begin{split}
\nabla \times \bf E+\partial_t \bf B&=0,\\
\nabla \times {\bf B}-\frac{1}{c^2}\partial_t \bf E&=\mu \bf J,\\
\nabla \cdot \bf E&=\frac{1}{\varepsilon}\rho,\\
\nabla \cdot \bf B&=0.
\end{split}
\end{equation}
We  now rewrite the Maxwell's equations into a form that is a suitable starting point for our EOS formulation of the elctromagnetic wave scattering problem.

First, observe that
\begin{equation}\label{ebrelation}
\begin{split}
&\partial_t \nabla \cdot {\bf B}=0,\\
&\partial_t \nabla \cdot {\bf E}=-\frac{1}{\varepsilon}\nabla \cdot {\bf J}.
\end{split}
\end{equation}
Equations (\ref{eq:maxwell}) and  (\ref{ebrelation}) lead to
\begin{equation}\label{erelation}
\partial_t(\nabla \cdot {\bf E}-\frac{1}{\varepsilon} \rho)=-\frac{1}{\varepsilon}(\partial_t \rho+\nabla \cdot {\bf J}).
\end{equation}
All fields we consider will be driven by the source that will
operate for some finite time interval. This means that at some time
in the past $t=t_0,$ we have
\begin{equation*}
\begin{split}
\nabla \cdot {\bf B}({\bf x},t_0)&=0,\\
\nabla \cdot {\bf E}({\bf x},t_0)&=\frac{1}{\varepsilon} \rho({\bf x},t_0)=0,
\end{split}
\end{equation*}
and this together with (\ref{ebrelation}) and (\ref{erelation}) imply that for any $t$ 
\begin{equation*}
\nabla \cdot {\bf B}({\bf x},t)=0
\end{equation*}
holds true. If we now use the equation of charge conservation $$
\partial_t\rho+\nabla \cdot {\bf J}=0
$$ then
\begin{equation*}
\nabla \cdot {\bf E}({\bf x},t)=\frac{1}{\varepsilon} \rho({\bf x},t)
\end{equation*} also holds true at all time. Taking these considerations into account, Maxwell's equations can be written in the following equivalent form
\begin{subequations}\label{premaxwell}
\begin{align}
&\nabla \times {\bf E}+\partial_t {\bf B}=0,\label{max1}\\
&\nabla \times {\bf B}-\frac{1}{c^2}\partial_t {\bf E}=\mu {\bf J},\label{max2}\\
&\partial_t\rho+\nabla \cdot {\bf J}=0.\label{max3}
\end{align}
\end{subequations}
 In order to complete the model, we must supply an equation of motion for the current $\bf J$
\begin{equation*}
\partial_t {\bf J}=F[{\bf J},\rho,{\bf E},{\bf B}].
\end{equation*}
The specific form for the functional $F$ is determined by what kind of material response we are considering. In order for the system to lead to an efficient numerical method it is important that the sources $\rho, {\bf J}$ are confined to some small region.  In this paper, in order to be specific, we look at the case of  a small metallic object interacting with light. We are not seeking to make a detailed computational investigation of this system, but is rather focused on testing  our computational approach with respect to implementational complexity and numerical stability. For this reason we choose the following simple nonlinear model for the metal response of such a system
\begin{equation}\label{Jmotion3}
\partial_t {\bf J}=(\alpha-\beta \rho) {\bf E}-\gamma {\bf J},
\end{equation}
where $\alpha, \beta$ and $\gamma$ are constants.

Following the usual approach, it is easy to show that the electric field satisfy the following equation

\begin{equation}\label{ewave}
\frac{1}{c^2}\partial_{tt} {\bf E}-\nabla^2{\bf E}=-\frac{1}{\varepsilon}\nabla\rho-\mu \partial_t{\bf J}.
\end{equation}
Each vector component of  equation (\ref{ewave}) is an inhomogeneous wave
equation. 
Let's suppose the scattering object is confined in a compact homogeneous domain denoted by $V_1$ while the light source is located in an unbounded domain outside the object which is denoted by $V_0.$ $\mu, \varepsilon$ are the magnetic permeability and the electric permittivity with their values $\mu_1, \varepsilon_1$ inside and $\mu_0, \varepsilon_0$ outside respectively.  $c$ represents the speed of light, with value $c_1$ inside  and $c_0$ outside the scattering object $V_1$. The sources ${\bf J_0}$ and $\rho_0$ are given and ${\bf J_1}$, $\rho_1$ are the response sources generated by the metallic object interacting with the light field.

We are now ready to start the construction of  the EOS formulations of this scattering problem. 

Applying the integral relation for the wave
equation (\ref{eq:mfundamental}) derived in Appendix \ref{integralidentity3d} on equation (\ref{ewave}) in domain $V_0$ and $V_1$ respectively, we get
\begin{equation}\label{eq:original e solution}
\begin{split}
{\bf E}_j({\bf x},t)&=-\int_{V_j}\,\mathrm{d}V'h_j({\bf x'},{\bf x})\{\mu_j \partial_{t'}{\bf J}_j+\frac{1}{\varepsilon_j}\nabla'\rho_j\}({\bf x'},T)\\
&\mp \int_S\,\mathrm{d}S'\{h_j({\bf x'},{\bf x})(\partial_{{\bf n}'}{\bf E}_j)({\bf x'},T)-\partial_{{\bf n}'}h_j({\bf x'},{\bf x}){\bf E}_j({\bf x'},T)\\
&+ \frac{1}{c_j}h_j({\bf x'},{\bf x})\partial_{{\bf n}'}|{\bf x'}-{\bf x}|(\partial_{t'}{\bf E}_j)({\bf x'},T)\},
\end{split}
\end{equation}
where $$
h_j({\bf x'},{\bf x})=\frac{c_j}{4\pi |{\bf x'}-{\bf x}|},
$$
with $j=0$ representing the outside domain $V_0$ and $j=1$ representing the inside domain $V_1$. Here ${\bf x}\in V_j$ and ${\bf n'}$ is the unit normal to the boundary, $S$ of $V_1$, at the point ${\bf x'}\in S$,   pointing out of the domain $V_1$.  The upper sign applies to the case $j=0$ and the lower sign for the case $j=1$. The same convention applies to all the following expressions in this section.

After a series of algebraic manipulations, starting with  (\ref{eq:original e solution}), we obtain
\begin{equation}\label{esolution}
\begin{split}
&{\bf E}_j({\bf x},t)=-\partial_t\frac{\mu_j}{4\pi}\int_{V_j}\,\mathrm{d}V'\frac{{\bf J}_j({\bf x'},T)}{|{\bf x'}-{\bf x}|}-\nabla\frac{1}{4\pi\varepsilon_j}\int_{V_j}\,\mathrm{d}V'\frac{\rho_j({\bf x'},T)}{|{\bf x'}-{\bf x}|}\\
&\mp \partial_t[\frac{1}{4\pi}\int_S\,\mathrm{d}S'\{\frac{1}{c_j|{\bf x'}-{\bf x}|} ({\bf n'}\times {\bf E}_j({\bf x'},T))\times \nabla'|{\bf x'}-{\bf x}|\\
&+\frac{1}{c_j|{\bf x'}-{\bf x}|} ({\bf n'}\cdot {\bf E}_j({\bf x'},T))\nabla'|{\bf x'}-{\bf x}|+\frac{1}{|{\bf x'}-{\bf x}|} {\bf n'}\times {\bf B}_j({\bf x'},T)  \}]\\
&\pm \frac{1}{4\pi}\int_S\,\mathrm{d}S'\{ ({\bf n'}\times {\bf E}_j({\bf x'},T))\times \nabla'\frac{1}{|{\bf x'}-{\bf x}|}\\
&+({\bf n'}\cdot {\bf E}_j({\bf x'},T))\nabla'\frac{1}{|{\bf x'}-{\bf x}|}\}.
\end{split}
\end{equation}
These manipulations are detailed  in Appendix \ref{Aintegral}.

\noindent  Like the electric field, the magnetic field also satisfies a wave equation
\begin{equation}\label{bwave}
\frac{1}{c^2}\partial_{tt} {\bf B}-\nabla^2{\bf B}=\mu \nabla \times {\bf J}.
\end{equation}
After a set of algebraic manipulations, similar to the ones we did for the electric field, we obtain
\begin{equation}\label{bsolution}
\begin{split}
&{\bf B}_j({\bf x},t)=\nabla\times\frac{\mu_j}{4\pi}\int_{V_j}\,\mathrm{d}V'\frac{{\bf J}_j({\bf x'},T)}{|{\bf x'}-{\bf x}|}\\
&+\partial_t[\frac{1}{4\pi}\int_S\,\mathrm{d}S'\{\frac{1}{c_j|{\bf x'}-{\bf x}|} ({\bf n'}\times {\bf B}_j({\bf x'},T))\times \nabla'|{\bf x'}-{\bf x}|\\
&\mp \frac{1}{c_j|{\bf x'}-{\bf x}|} ({\bf n'}\cdot {\bf B}_j({\bf x'},T))\nabla'|{\bf x'}-{\bf x}|-\frac{1}{c_j^2}\frac{1}{|{\bf x'}-{\bf x}|} {\bf n'}\times {\bf E}_j({\bf x'},T)  \}]\\
&\pm \frac{1}{4\pi}\int_S\,\mathrm{d}S'\{ ({\bf n'}\times {\bf B}_j({\bf} x',T))\times \nabla'\frac{1}{|{\bf x'}-{\bf x}|}\\
&+({\bf n'}\cdot {\bf B}_j({\bf x'},T))\nabla'\frac{1}{|{\bf x'}-{\bf x}|}\}.
\end{split}
\end{equation}
The identities (\ref{esolution}) and (\ref{bsolution}), for the electric and magnetic field, are our version of the general integral identities for the electromagnetic field derived by D.S.Jones\cite{Jones64}.
In addition to these two identities we get, in a very similar way,  two additional integral identities \cite{Jakobsen09},
\begin{equation}\label{e0solution}
\begin{split}
&0=-\partial_t\frac{\mu_{1-j}}{4\pi}\int_{V_{1-j}}\,\mathrm{d}V'\frac{{\bf J}_{1-j}({\bf x'},T)}{|{\bf x'}-{\bf x}|}-\nabla\frac{1}{4\pi\varepsilon_{1-j}}\int_{V_{1-j}}\,\mathrm{d}V'\frac{\rho_{1-j}({\bf x'},T)}{|{\bf x'}-{\bf x}|}\\
&\pm \partial_t[\frac{1}{4\pi}\int_S\,\mathrm{d}S'\{\frac{1}{c_{1-j}|{\bf x'}-{\bf x}|} ({\bf n'}\times {\bf E}_{1-j}({\bf x'},T))\times \nabla'|{\bf x'}-{\bf x}|\\
&+\frac{1}{c_{1-j}|{\bf x'}-{\bf x}|} ({\bf n'}\cdot {\bf E}_{1-j}({\bf x'},T))\nabla'|{\bf x'}-{\bf x}|+\frac{1}{|{\bf x'}-{\bf x}|} {\bf n'}\times {\bf B}_{1-j}({\bf x'},T)  \}]\\
&\mp \frac{1}{4\pi}\int_S\,\mathrm{d}S'\{ ({\bf n'}\times {\bf E}_{1-j}({\bf x'},T))\times \nabla'\frac{1}{|{\bf x'}-{\bf x}|}\\
&+({\bf n'}\cdot {\bf E}_{1-j}({\bf x'},T))\nabla'\frac{1}{|{\bf x'}-{\bf x}|}\},
\end{split}
\end{equation}
and
\begin{equation}\label{b0solution}
\begin{split}
&0=\nabla\times\frac{\mu_{1-j}}{4\pi}\int_{V_{1-j}}\,\mathrm{d}V'\frac{{\bf J}_{1-j}({\bf x'},T)}{|{\bf x'}-{\bf x}|}\\
&\pm \partial_t[\frac{1}{4\pi}\int_S\,\mathrm{d}S'\{\frac{1}{c_{1-j}|{\bf x'}-{\bf x}|} ({\bf n'}\times {\bf B}_{1-j}({\bf x'},T))\times \nabla'|{\bf x'}-{\bf x}|\\
&+\frac{1}{c_{1-j}|{\bf x'}-{\bf x}|} ({\bf n'}\cdot {\bf B}_{1-j}({\bf x'},T))\nabla'|{\bf x'}-{\bf x}|-\frac{1}{c_{1-j}^2}\frac{1}{|{\bf x'}-{\bf x}|} {\bf n'}\times {\bf E}_{1-j}({\bf x'},T)  \}]\\
&\mp \frac{1}{4\pi}\int_S\,\mathrm{d}S'\{ ({\bf n'}\times {\bf B}_{1-j}({\bf} x',T))\times \nabla'\frac{1}{|{\bf x'}-{\bf x}|}\\
&+({\bf n'}\cdot {\bf B}_{1-j}({\bf x'},T))\nabla'\frac{1}{|{\bf x'}-{\bf x}|}\},
\end{split}
\end{equation}
for ${\bf x} \in V_j,\  j=0, 1.$
In the above expressions, $${\bf E}_j({\bf x}', t)=\lim_{{\bf x}\rightarrow {\bf x}'}{\bf E}_j({\bf x}, t)$$ and 
 $${\bf B}_j({\bf x}', t)=\lim_{{\bf x}\rightarrow {\bf x}'}{\bf B}_j({\bf x}, t),$$
where ${\bf x} \in V_j, \, j=0, 1.$
In the end, we have a full set of the integral identities of  the inside and the outside fields expressed by (\ref{esolution}), (\ref{bsolution}), (\ref{e0solution}) and (\ref{b0solution}) which can be written compactly as
\begin{equation}\label{compactv}
\begin{split}
{\bf E}_1&=M_1({\bf n}'\times {\bf E}_1,{\bf n}'\cdot {\bf E}_1,{\bf n}'\times {\bf B}_1),\\
0&=M_0({\bf n}'\times {\bf E}_0,{\bf n}'\cdot {\bf E}_0,{\bf n}'\times {\bf B}_0),\\
{\bf B}_1&=N_1({\bf n}'\times {\bf B}_1,{\bf n}'\cdot {\bf B}_1,{\bf n}'\times {\bf E}_1),\\
0&=N_0({\bf n}'\times {\bf B}_0,{\bf n}'\cdot {\bf B}_0,{\bf n}'\times {\bf E}_0),
\end{split}
\end{equation}
for ${\bf x} \in V_1$ and 
\begin{equation}\label{compactv0}
\begin{split}
{\bf E}_0&=M_0({\bf n}'\times {\bf E}_0,{\bf n}'\cdot {\bf E}_0,{\bf n}'\times {\bf B}_0),\\
0&=M_1({\bf n}'\times {\bf E}_1,{\bf n}'\cdot {\bf E}_1,{\bf n}'\times {\bf B}_1),\\
{\bf B}_0&=N_0({\bf n}'\times {\bf B}_0,{\bf n}'\cdot {\bf B}_0,{\bf n}'\times {\bf E}_0),\\
0&=N_1({\bf n}'\times {\bf B}_1,{\bf n}'\cdot {\bf B}_1,{\bf n}'\times {\bf E}_1,
\end{split}
\end{equation}
for ${\bf x} \in V_0.$

We will now derive the boundary integral identities of (\ref{premaxwell}) by letting ${\bf x}$ approach  the surface from the inside and the outside of the scattering object $V_1$, separately. We observe that, in this limit, strong singularities only appear in the last term of the integrals  in (\ref{compactv}) and (\ref{compactv0}). Hence we are faced with a singular term which takes the form of
\begin{equation}\label{Iequation}
\begin{split}
{\bf I}&= \lim_{\epsilon \rightarrow 0}\int_{S_\epsilon}\,\mathrm{d}S'\{ ({\bf n'}\times {\bf A}({\bf x'},T))\times \nabla'\frac{1}{|{\bf x'}-{\bf x}|}+({\bf n'}\cdot {\bf A}({\bf x'},T))\nabla'\frac{1}{|{\bf x'}-{\bf x}|}\}\\
&=\lim_{\epsilon\rightarrow 0}\int_{S_\epsilon}\,\mathrm{d}S'\{ (\frac{{\bf x'}-{\bf x}}{|{\bf x'}-{\bf x}|^3}\cdot {\bf A}){\bf n}'-(\frac{{\bf x'}-{\bf x}}{|{\bf x'}-{\bf x}|^3}\cdot{\bf n}' ){\bf A} -({\bf n}'\cdot {\bf A})\frac{{\bf x'}-{\bf x}}{|{\bf x'}-{\bf x}|^3}\},
\end{split}
\end{equation}
where ${\bf A}({\bf x'},T)$ is a vector function with ${\bf x'}$  located on the surface $S_\epsilon$ which is an small disk of radius $\epsilon.$
If we let ${\bf x}$ approach  a surface point ${\boldsymbol \xi}$, from the inside of $V_1$, along a direction
$$ {\bf x}-{\boldsymbol \xi}=\epsilon {\bf a}=-\epsilon \alpha {\bf n}-\epsilon{\boldsymbol \beta}, $$
where $ {\bf n}$ is the unit normal vector pointing out of $V_1$, at the point ${\boldsymbol \xi}$, and ${\boldsymbol \beta}$ is a unit vector along the direction ${\bf x}'-{\boldsymbol \xi}$, tangential to $S$,  at the same point ${\boldsymbol \xi},$ 
we have 
\begin{equation}\label{singularrelation}
\begin{split}
\lim_{\epsilon\rightarrow 0}\int_{S_\epsilon} \frac{{\bf x'}-{\bf x}}{|{\bf x'}-{\bf x}|^3}\,\mathrm{d}S=\lim_{\epsilon \rightarrow 0}\int_{S_\epsilon}\frac{\boldsymbol \eta+\epsilon \alpha {\bf n}}{|\boldsymbol \eta+\epsilon \alpha {\bf n}|^3} \,\mathrm{d}S,
\end{split}
\end{equation}
where $\eta={\bf x}'-{\boldsymbol \xi}+\epsilon{\boldsymbol \beta}.$
Using spherical coordinates, (\ref{singularrelation}) turns into
\begin{equation*}
\begin{split}
\lim_{ \epsilon \rightarrow 0}\int_0^{2\pi}\int_0^{\epsilon}\rho \frac{(\rho \cos\theta,\rho\sin\theta,\epsilon\alpha)}{({\rho}^2+(\epsilon\alpha)^2)^{\frac{3}{2}}}\,\mathrm{d}\theta\,\mathrm{d}\rho=\chi{\bf n},
\end{split}
\end{equation*}
where $\chi=2\pi\alpha(1-\frac{1}{\sqrt{\alpha^2+1}})$.
Similarly, if ${\bf x}$ approaches  ${\boldsymbol \xi}$ from outside of $V_1$, we have,
$$ \lim_{ \epsilon \rightarrow 0}\int_{S_\epsilon} \frac{{\bf x'}-{\bf x}}{|{\bf x'}-{\bf x}|^3}\,\mathrm{d}S=-\chi{\bf n}.$$
So in the end, 
\begin{equation*}
\begin{split}
I_+&=\chi{\bf A},\\
I_-&=-\chi{\bf A},
\end{split}
\end{equation*}
where $I_+$ and $I_-$ denote the limit of equation (\ref{Iequation}) by letting ${\bf x}$ approach ${S_\epsilon} $ from the inside and the outside of $V_1$ respectively.
After taking these inside and outside limits, we get the following set of equations 
\begin{equation}\label{eset}
\begin{split}
{\bf E}_+&=M_1({\bf n}'\times {\bf E}_+,{\bf n}'\cdot {\bf E}_+,{\bf n}'\times {\bf B}_+)+\chi{\bf E}_+,\\
0&=M_0({\bf n}'\times {\bf E}_-,{\bf n}'\cdot {\bf E}_-,{\bf n}'\times {\bf B}_-)-\chi{\bf E}_-,\\
{\bf E}_-&=M_0({\bf n}'\times {\bf E}_-,{\bf n}'\cdot {\bf E}_-,{\bf n}'\times {\bf B}_-)+\chi{\bf E}_-,\\
0&=M_1({\bf n}'\times {\bf E}_+,{\bf n}'\cdot {\bf E}_+,{\bf n}'\times {\bf B}_+)-\chi{\bf E}_+,
\end{split}
\end{equation}
where ${\bf E}_+$ is the limit of ${\bf E_1}$ with ${\bf x}$ approaching the surface from the inside of the object while ${\bf E}_-$ is the limit of ${\bf E_0}$ with ${\bf x}$ approaching the surface from the outside of the object. These equations, because of the limits taken, contains singular integrals that must be interpreted as Cauchy principal integrals.
Adding the first two equations of (\ref{eset}) gives us
\begin{equation}\label{eplus}
\begin{split}
{\bf E}_+&=M_1({\bf n}'\times {\bf E}_+,{\bf n}'\cdot {\bf E}_+,{\bf n}'\times {\bf B}_+)+M_0({\bf n}'\times {\bf E}_-,{\bf n}'\cdot {\bf E}_-,{\bf n}'\times {\bf B}_-)\\
&+\chi{\bf E}_+-\chi{\bf E}_-,
\end{split}
\end{equation}
and adding the last two equations of (\ref{eset}) gives us 
\begin{equation}\label{eminus}
\begin{split}
{\bf E}_-&=M_0({\bf n}'\times {\bf E}_-,{\bf n}'\cdot {\bf E}_-,{\bf n}'\times {\bf B}_-)+M_1({\bf n}'\times {\bf E}_+,{\bf n}'\cdot {\bf E}_+,{\bf n}'\times {\bf B}_+)\\
&+\chi{\bf E}_--\chi{\bf E}_+.
\end{split}
\end{equation}
Repeating the derivations we just did for the electric field, give us, in a similar way, the following set of equations for the magnetic field 
\begin{equation}\label{bplus}
\begin{split}
{\bf B}_+&=N_1({\bf n}'\times {\bf B}_+,{\bf n}'\cdot {\bf B}_+,{\bf n}'\times {\bf E}_+)+N_0({\bf n}'\times {\bf B}_-,{\bf n}'\cdot {\bf B}_-,{\bf n}'\times {\bf E}_-)\\
&+\chi{\bf B}_+-\chi{\bf B}_-,
\end{split}
\end{equation}
\begin{equation}\label{bminus}
\begin{split}
{\bf B}_-&=N_0({\bf n}'\times {\bf B}_-,{\bf n}'\cdot {\bf B}_-,{\bf n}'\times {\bf E}_-)+N_1({\bf n}'\times {\bf B}_+,{\bf n}'\cdot {\bf B}_+,{\bf n}'\times {\bf E}_+)\\
&+\chi{\bf B}_--\chi{\bf B}_+,
\end{split}
\end{equation}
where ${\bf B}_+$ is the limit of ${\bf B_1}$ with ${\bf x}$ approaching the surface from the inside of the object while ${\bf B}_-$ is the limit of ${\bf B_0}$ with ${\bf x}$ approaching the surface from the outside of the object. Also in these equations the singular integrals that occur must be interpreted as Cauchy principal value integrals.
So far,  we have two outer equations for the outer limit fields ${\bf E}_-$, ${\bf B}_-$ and two inner equations for the inner limit fields ${\bf E}_+$, ${\bf B}_+$. We also have the usual electromagnetic boundary conditions at the surface $S$ which separate regions with different susceptibilities and permittivities
\begin{equation*}
\begin{split}
{\bf n}'\times {\bf E}_+&={\bf n}'\times {\bf E}_-,\\
{\bf n}'\times {\bf B}_+&=\frac{u_1}{u_0}{\bf n}'\times {\bf B}_-,\\
{\bf n}'\cdot {\bf B}_+&={\bf n}'\cdot {\bf B}_-,\\
{\bf n}'\cdot {\bf E}_+&=\frac{\varepsilon_0}{\varepsilon_1}{\bf n}'\cdot {\bf E}_-.
\end{split}
\end{equation*}
It might appear that we have more equations than we need here. The very same problem was encountered earlier while deriving the EOS formulation for two 1D toy models \cite{Aihua}.  It appears as if we can
use the two outer equations to solve for ${\bf E}_-$ and
${\bf B}_-$ and then use the boundary conditions to find
${\bf E}_+$ and ${\bf B}_+$. But these field values inside the scattering object cannot in general be consistent with the field values derived directly from the two inner equations for  ${\bf E}_+$ and ${\bf B }_+$. For example, if there is a source  inside of $V_1$ and no source outside of $V_1$,  the first approach would give vanishing electric and magnetic field whereas the second approach certainly would not. On the other hand, for a given source, the Maxwell equations has a unique solution, which by construction also satisfy all the integral identities.

In order to understand what the problem is, and how to fix it, we must just realize that, from an abstract point of view, we have the following formal situation
\begin{equation}\label{eq:singular system}
\begin{split}
A{\bf X} &={\bf b},\\
B{\bf X} &={\bf c},
\end{split}
\end{equation}
where A and  B are singular but where we know that  (\ref{eq:singular system}) has a unique solution. 

 In this situation, let us assume that $\alpha A+ B$ is nonsingular for some choice of $\alpha.$ Any solution of (\ref{eq:singular system}) is a solution of 
\begin{equation}\label{eq:nonsingular equation}
(\alpha A+B){\bf X}=\alpha {\bf b}+{\bf c},
\end{equation}
and since $\alpha A+ B$ is nonsingular the unique solution of the singular system (\ref{eq:singular system}) must in fact be the unique solution of the nonsingular system (\ref{eq:nonsingular equation}). We know that the
solution of Maxwell is unique for a given source, so since the
integral equations are equivalent to Maxwell, our four integral
equations for the two unknown fields on $S$  must have a unique
solution. This happens only if they are singular. Thus in our situation, we can simply add (\ref{eplus}) and (\ref{eminus}), which gives
\begin{equation}\label{epm}
\begin{split}
{\bf E}_++{\bf E}_-&=2(M_1({\bf n}'\times {\bf E}_+,{\bf n}'\cdot {\bf E}_+,{\bf n}'\times {\bf B}_+)\\
&+M_0({\bf n}'\times {\bf E}_-,{\bf n}'\cdot {\bf E}_-,{\bf n}'\times {\bf B}_-)),
\end{split}
\end{equation}
and add (\ref{bplus}) and (\ref{bminus}), which gives
\begin{equation}\label{bpm}
\begin{split}
{\bf B}_++{\bf B}_-&=2(N_1({\bf n}'\times {\bf B}_+,{\bf n}'\cdot {\bf B}_+,{\bf n}'\times {\bf E}_+)\\
&+N_0({\bf n}'\times {\bf B}_-,{\bf n}'\cdot {\bf B}_-,{\bf n}'\times {\bf E}_-)).
\end{split}
\end{equation}

Observe that for any vector ${\bf A}$, the following identities hold true
\begin{equation}\label{id1}
\begin{split}
&{\bf n}\times({\bf n}\times {\bf A})=({\bf n}\cdot {\bf A}){\bf n}-({\bf n}\cdot {\bf n}){\bf A},\\
&{\bf A}=({\bf n}\cdot {\bf A}){\bf n}-{\bf n}\times({\bf n}\times {\bf A}).
\end{split}
\end{equation}
Performing (\ref{id1}) on (\ref{epm}) and (\ref{bpm}) we obtain the following final boundary integral identities
\begin{subequations}\label{maxboundary}
\begin{align}
(I+\frac{1}{2}(\frac{\varepsilon_1}{\varepsilon_0}-1){\bf n}\ {\bf n}){\bf E}_+({\bf x},t)&={\bf I}_e+{\bf O}_e+{\bf B}_e,\label{eb}\\
(I+\frac{1}{2}(1-\frac{\mu_0}{\mu_1}){\bf n}\ {\bf n}){\bf B}_+({\bf x},t)&={\bf I}_b+{\bf O}_b+{\bf B}_b,\label{bb}
\end{align}
\end{subequations}
where ${\bf x} \in S,$  ${\bf n}$ is the unit normal vector pointing out of $V_1$ at the point ${\bf x},$ $I$ is the identity matrix and
\begin{equation}
{\bf I}_e=-\partial_t\frac{\mu_1}{4\pi}\int_{V_1}\,\mathrm{d}V'\frac{{\bf J}_1({\bf x'},T)}{|{\bf x'}-{\bf x}|}-\frac{1}{4\pi\varepsilon_1}\nabla\int_{V_1}\,\mathrm{d}V'\frac{\rho_1({\bf x'},T)}{|{\bf x'}-{\bf x}|},\label{ie}
\end{equation}
\begin{equation}
{\bf O}_e=-\partial_t\frac{\mu_0}{4\pi}\int_{V_0}\,\mathrm{d}V'\frac{{\bf J}_0({\bf x'},T)}{|{\bf x'}-{\bf x}|}-\nabla\frac{1}{4\pi\varepsilon_0}\int_{V_0}\,\mathrm{d}V'\frac{\rho_0({\bf x'},T)}{|{\bf x'}-{\bf x}|},\label{oe}\\
\end{equation}
\begin{equation}\label{Be}
\begin{split}
&{\bf B}_e=\frac{1}{4\pi}\partial_t\int_{S}\,\mathrm{d}S^{'}\{(\frac{1}{c_1}-\frac{1}{c_0})\frac{1}{|{\bf x'}-{\bf x}|} ({\bf n}^{'}\times {\bf E}_+({\bf x'},T))\times \nabla'|{\bf x'}-{\bf x}|\\
&\indent +(\frac{1}{c_1}-\frac{\varepsilon_1}{\varepsilon_0 c_0})\frac{1}{|{\bf x'}-{\bf x}|} ({\bf n}^{'}\cdot {\bf E}_+({\bf x'},T))\nabla'|{\bf x'}-{\bf x}|\\
&\indent+(1-\frac{\mu_0}{\mu_1})\frac{1}{|{\bf x'}-{\bf x}|} {\bf n}^{'}\times {\bf B}_+({\bf x'},T)\} \\
&\indent-\frac{1}{4\pi}\int_{S}\,\mathrm{d}S^{'}((1-\frac{\varepsilon_1}{\varepsilon_0})({\bf n}^{'}\cdot {\bf E}_+({\bf x'},T))\nabla'\frac{1}{|{\bf x'}-{\bf x}|}),
\end{split}
\end{equation}
\begin{equation}
{\bf I}_b=\nabla\times\frac{\mu_1}{4\pi}\int_{V_1}\,\mathrm{d}V'\frac{{\bf J}_1({\bf x'},T)}{|{\bf x'}-{\bf x}|},\label{ib}\\
\end{equation}
\begin{equation}
{\bf O}_b=\nabla\times\frac{\mu_0}{4\pi}\int_{V_0}\,\mathrm{d}V'\frac{{\bf J}_0({\bf x'},T)}{|{\bf x'}-{\bf x}|},\label{ob}
\end{equation}
\begin{equation}\label{Bb}
\begin{split}
&{\bf B}_b=\frac{1}{4\pi}\partial_t\int_{S}\,\mathrm{d}S^{'}\{(\frac{1}{c_1}-\frac{\mu_0}{\mu_1 c_0})\frac{1}{|{\bf x'}-{\bf x}|} ({\bf n}^{'}\times {\bf B}_+({\bf x'},T))\times \nabla'|{\bf x'}-{\bf x}|\\
&\indent+(\frac{1}{c_1}-\frac{1}{c_0})\frac{1}{|{\bf x'}-{\bf x}|} ({\bf n}^{'}\cdot {\bf B}_+({\bf x'},T))\nabla'|{\bf x'}-{\bf x}|\\
&\indent+(\frac{1}{c_0^2}-\frac{1}{c_1^2})\frac{1}{|{\bf x'}-{\bf x}|} {\bf n}^{'}\times {\bf E}_+({\bf x'},T)\} \\
&\indent-\frac{1}{4\pi}\int_{S}\,\mathrm{d}S^{'}\{(1-\frac{u_0}{u_1}) ({\bf n}^{'}\times {\bf B}_+({\bf x'},T))\times \nabla'\frac{1}{|{\bf x'}-{\bf x}|}\}.
\end{split}
\end{equation}
Note that ${\bf O}_e$ and ${\bf O}_b$ are fields on the surfaces generated by the source in the absence of the scattering objects.
Equations (\ref{premaxwell})  and (\ref{Jmotion3}) together with the boundary integral identities (\ref{eb}) and (\ref{bb}) compose the EOS formulations  of our model.
\section{Artificial source test and numerical implementation}\label{implementation3d}
The motivation, for introducing the EOS formulation for Maxwell's equations, is a numerical one. The technical issues occuring for the numerical implementation discussed  in this paper, will be part of any numerical implementation of our scheme, which in general will involve multiple, arbitrarily shaped, scattering objects,  that include linear and nonlinear optical response.  We expect however, that the general nature of these issues will reveal themselves already in the simplest possible setting, where we have one scattering object of rectangular shape.
The numerical implementation consists of a domain method for the model (\ref{premaxwell})  and (\ref{Jmotion3}),  determining the evolution of the fields inside the scattering object, and a scheme for updating the boundary values of the fields using the integral identities (\ref{eb}) and (\ref{bb}).
For the internal domain method we choose to use a combination of  Lax-Wendroff on (\ref{premaxwell}) and modified Euler's method on  (\ref{Jmotion3}), this is similar to what we did for the simple 1D case \cite{Aihua},  previously. For the boundary part of the scheme we use the mid-point rule to the nonsingular integrals appearing in (\ref{maxboundary}). The treatment of the singular integrals is technical and rather lengthy and will therefore be reported elsewhere  \cite{Aihua3}. Here it is enough to note that we calculate the singular integrals by reducing them to a singular core, which we calculate exactly, and nonsingular surface and line integrals, that we calculate numerically. The reductions proceed through a nontrivial  use of well known integral identities.

For the inside of the scattering object we will use a rectangular grid. This grid is however not uniform close to the boundary. This is because the grid has to support both the discrete approximations to the partial derivatives and discrete approximations  to the integrals, used to update the boundary values based on the current and previous internal values of the fields. The fact that in our scheme the boundary values are dynamical variables, enforce some special difference rules that applies close to the boundary. This is an extra complication for our scheme, but they are manageable, and will be part of any scheme that implements the EOS formulation introduced in this paper. Details are given in \cite{Aihua3}.

What we do in this section is to report on some tests that we have run on our scheme. A usual approach to testing of numerical implementations involve finding exact special solutions corresponding to special source functions.  In this section we do not use this approach, but rather use an artificial source test to verify the correctness of our EOS formulations. The basic idea behind the artificial source test, of some numerical scheme designed for a system of PDEs, $\mathscr{L}\psi=0$ ,  is to slightly modify the system by adding an arbitrary source to all the equations in the system, creating a new modified system $\mathscr{L}\psi=g$ . This modification typically lead to minimal modifications to the numerical scheme, where most of the effort and complexity are usually spent on the derivatives and nonlinear terms. For the equations, however, the presence of the sources  change the situation completely. This is because the presence of the added sources implies that {\it any} function is a solution to the equations for {\it some} choice of sources.
Thus we can pick a function $\psi_0$ and insert it into the model and  calculate the source function $g_0=\mathscr{L}\psi_0$ so that our chosen function is a solution to the extended equation. Finally we run the numerical scheme with the calculated source function and compare the numerical solution with the exact solution $\psi_0$.

\begin{figure}[h!]
\centering
\includegraphics[width=10cm]{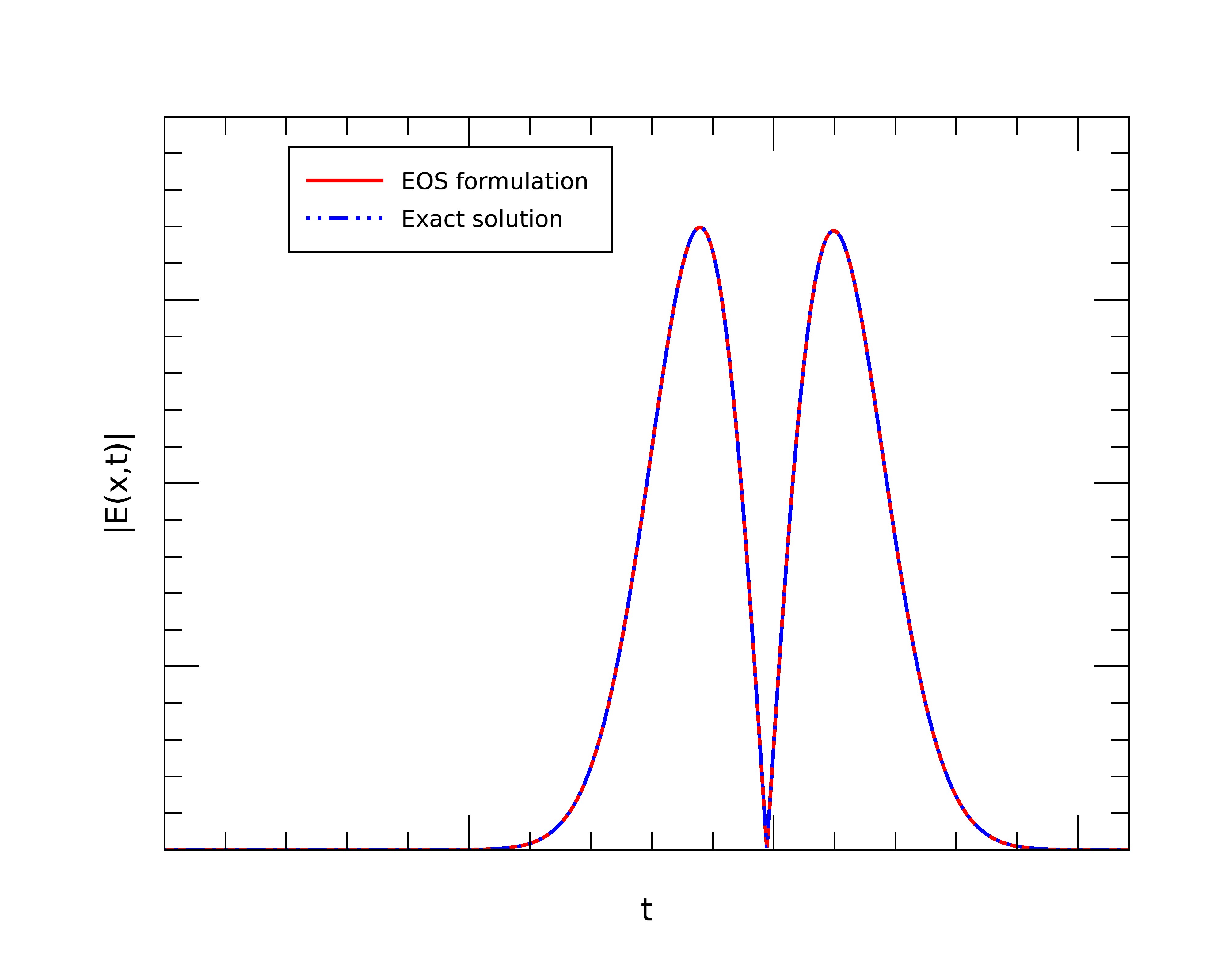}
\caption{Artificial source test: $ |{\bf E}({\bf x}, t)|.$  Intensity of electric field at a specific point at different times.  $ b=1.0, \alpha_1=40, \beta_1=1.0, x_o=0.0, y_o=0.0, z_o=0.0, t_a=1.0, c_0=1.0,\mu_0=1.0, \varepsilon_0=1.0, c_1=0.82, \mu_1=1.0, \varepsilon_1=1.5, \tau=0.45, \alpha=1.0, \beta=0.01, \gamma=0.01.$}
\label{3dtest}
\end{figure}

A modified model of (\ref{premaxwell}) and (\ref{Jmotion3}) with artificial sources is generally constructed by
\begin{equation*}
\begin{split}
&\partial_t {\bf B}+\nabla \times {\bf E}={\boldsymbol \varphi}_1,\\\nonumber
&\frac{1}{c_1^2}\partial_t {\bf E}-\nabla \times {\bf B}=-\mu_1 {\bf J}+{\boldsymbol \varphi}_2,\\\nonumber
&\nabla \cdot {\bf E}=\frac{1}{\varepsilon_1}\rho +{\varphi}_3,\\\nonumber
&\nabla \cdot  {\bf B}={\boldsymbol \varphi}_4,\\\nonumber
&\partial_t {\bf J}=(\alpha-\beta \rho)  {\bf
E}-\gamma {\bf J}+{{\boldsymbol \varphi}_5},\nonumber
\end{split}
\end{equation*}
where ${\boldsymbol \varphi}_1, {\boldsymbol \varphi}_2, {\boldsymbol \varphi}_4 $ and ${\boldsymbol \varphi}_5$ are a set of vector functions and $ {\varphi}_3$ is a scale function.
Observe that
$$\nabla \cdot (\nabla \times {\bf E}) +\nabla \cdot \partial_t {\bf B}=\nabla \cdot  {{\boldsymbol \varphi}_1},$$
and this yields $$\partial_t {\boldsymbol \varphi}_4=\nabla \cdot  {{\boldsymbol \varphi}_1}.$$
Based on this we suppose $ {{\boldsymbol \varphi}_1}=0$ and $ {{\boldsymbol \varphi}_4}=0$ which can simplify the choice of the exact solutions. We also observe that if $\varphi_2 $ and $\varphi_3$ are set to be both 0, then the continuity equation 
\begin{equation*}
\partial_t\rho+\nabla\cdot{\bf J}=0
\end{equation*}
 is automatically satisfied. 
So in the end, the source extended model is given by 
\begin{subequations}\label{amaxwell}
\begin{align}
&\partial_t {\bf B}+\nabla \times {\bf E}=0,\label{am1}\\
&\frac{1}{c_1^2}\partial_t {\bf E}-\nabla \times {\bf B}=-\mu_1 {\bf J},\label{am2}\\
&\nabla \cdot {\bf E}=\frac{1}{\varepsilon_1}\rho,\label{am3}\\
&\partial_t {\bf J}=(\alpha-\beta \rho)  {\bf
E}-\gamma {\bf J}+{{\boldsymbol \varphi}}.\label{am4}
\end{align}
\end{subequations}
For model (\ref{amaxwell}), any choice of $ \tilde {\bf E},  \tilde{\bf B} $ can be a solution if the artificial source is given by 
$$
{{\boldsymbol \varphi}}=\partial_t \tilde {\bf J}-(\alpha-\beta \tilde \rho)  \tilde {\bf E}-\gamma \tilde {\bf J},
$$
where $\tilde{\bf J}$ and $\tilde{\rho}$ are given respectively by 
\begin{equation*}
\begin{split}
& \tilde{{\bf J}}=\frac{1}{\mu_1}(\nabla \times \tilde{{\bf B}}-\frac{1}{c_1^2}\partial_t \tilde{{\bf E}}),\\
&\tilde{\rho}=\varepsilon_1 \nabla \cdot \tilde{{\bf E}}.
\end{split}
\end{equation*}
Due to (\ref{am1}), we can simply choose a vector function ${\boldsymbol \phi}$, such that 
\begin{equation}\label{exactsolution}
\begin{split}
\tilde{{\bf E}}&=-\partial_t {{\boldsymbol \phi}},\\
\tilde{{\bf B}}&=\nabla \times {\boldsymbol \phi}.
\end{split}
\end{equation}
Figure \ref{3dtest} shows the comparison between the numerical implementations and the exact solutions where we have used
\begin{equation*}
\boldsymbol \phi=(\arctan(b^2t^2)e^{-\alpha_1(x-x_o+y-y_o+z-z_o+\beta_1(t-t_a))^2}, 0, 0).\\
\end{equation*}
Values of the parameters are listed below the figure. From figure  it is evident that the agreement between the exact solution and the numerical solution is excellent.

For a general case where the electromagnetic fields inside the object are produced by the outside source,
we set up the outside source ${\bf J}_0$ and $\rho_0$ to be a combination of a bump function in time and a delta function in space which is easily integrated in space.
In order to satisfy the  continuity equation
$$
\partial_t \rho_0+\nabla \cdot {\bf J}_0=0,
$$
we can choose a vector function $\boldsymbol \varphi$  such that
\begin{equation}\label{j0p0}
\begin{split}
 {\bf J}_0&=-\partial_t {\boldsymbol \varphi},\\
\rho_0&=\nabla \cdot {\boldsymbol \varphi}.
\end{split}
\end{equation}
\begin{figure}[h!]
\centering
\includegraphics[width=10cm]{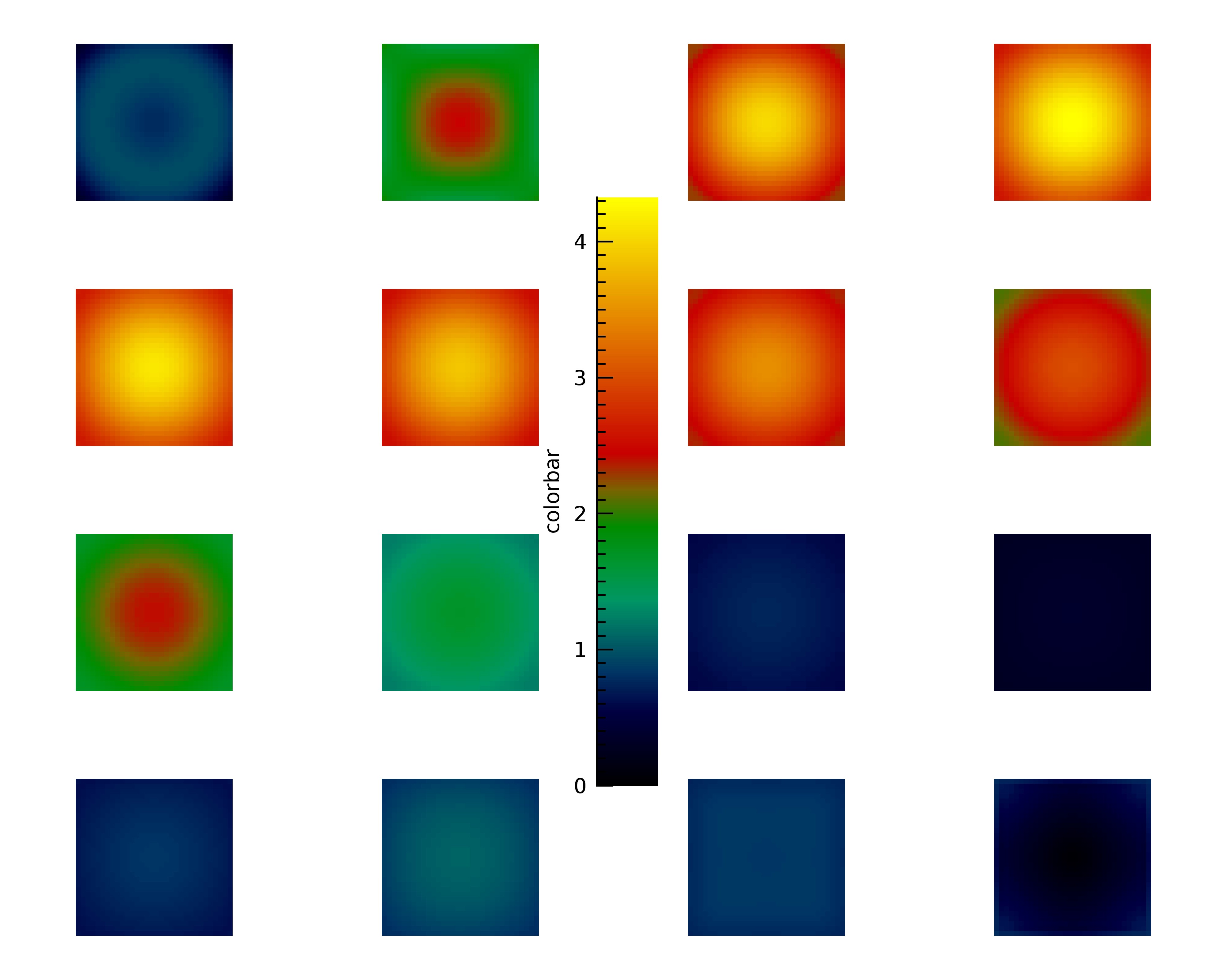}
\caption{The intensity of the electric field on a specific surface in $yz$ plane at different times. $t_0=1.5,  x_0=-0.3,  y_0=0.0,  z_0=0.0, c_0=1.0, \mu_0=1.0,  \varepsilon_0=1.0,  c_1=0.82,  \mu_1=1.0,  \varepsilon_1=1.5,  \tau=0.45, \alpha=0.1, \beta=0.01, \gamma=2.$ }
\label{3doriginal}
\end{figure}

Figure \ref{3doriginal} shows the intensity of the electric field on a surface in $yz$ plane at different times. Values of parameters used are shown below the figure. The figure shows a pulse of light passing through the plane, which is what we would expect from the nature of our chosen source.
For figure  \ref{3doriginal} , we have chosen 
$$ {\boldsymbol \varphi}=(\varphi, 0, 0),$$
where 
\[ \varphi ({\bf x}, t)=\begin{cases} 
      \delta({\bf x}-{\bf x}_0) e^{\frac{1}{(t-t_0)^2-1}}& t \in [t_0-1,t_0+1], \\
      0 & t \notin [t_0-1,t_0+1],
   \end{cases}
\]
with ${\bf x}_0=(x_0, y_0, z_0).$ 

Our numerical scheme, being explicit, is not unconditionally stable. There  is however a stable range,  $\tau_1< \tau <\tau_2$ , for the time step $\Delta t$  $$\Delta t=\frac{\tau}{c_1} Min \{\Delta x, \Delta y, \Delta z\},$$ where $\tau_1$ and $\tau_2$  determine the lower and upper boundaries of the stability range for the sceheme. We have carefully investigated the source of the upper and lower bound of the range and how the width of the stability range depends on material parameters. It is not appropriate to include these fairly techincal numerical investigations here, a full discussion will be presented elsewhere in \cite{Aihua3}. Here it is enough to note that the source of the lower stability bound is the numerical implementation of the boundary part of the algoritm, and the source of the upper bound is the numerical implementation of the domain part of the algorithm, for our case this is a combination of Lax-Wendroff for the electromagnetic fields, and modified Euler for the current .


\section{Conclusions}\label{conclusion3d}
In this paper we have showed that our EOS formulation of electromagnetic scattering can be accurately and stably implemented using one particular choice of numerical scheme for the inside of the objects and for the integral representation of the boundary values required by the inside scheme. For a stable numerical solution, the time step needs to be confined in some range, where we have found that this range is not only  determined by the internal domain-base method due to the non-uniform grids but also determined by the boundary integral representations. Discussions on how the internal non-uniform grids and the boundary integral representation effect the time stable range is reported in \cite{Aihua} which will be published elsewhere.
It is worth stressing that the existence of the stability range and its width depends not only on the material parameters but certainly on choices made for the numerical implementation of the boundary part and domain part of the algorithm. In principle, any numerical scheme can be used for the domain part of the algorithm, also the extremely well established FDTD method. It would be interesting to see how this method would perform with respect to stability. There is also the question of going fully implicit, both for the boundary part and the domain part of the algorithm. One would think that this would have a chance of producing an unconditionally stable algorithm for our EOS formulation for Maxwell's equations.
\begin{appendices}

\section{The integral identity for a 3D wave equation}\label{integralidentity3d}
We will start by considering a wave equation in 3D
\begin{equation}\label{3dwaveequation}
\frac{1}{c^2}\partial_{tt}\varphi({\bf x},t)-\nabla^2\varphi({\bf x},t)=\rho({\bf x},t),
\end{equation}
where ${\bf x}=(x,y,z),$ $c$ is the propagation speed and $$\nabla^2=\partial_x^2+\partial_y^2+\partial_z^2.$$ Let $D\times T$ be a given space-time domain. 
We will assume that the source $\rho({\bf x},t)$ is entirely contained in $D \times T$. The operator $$\mathcal{L}=\frac{1}{c^2}\partial_{tt}-\nabla^2$$ is formally self adjoint. Observe that for any pair of functions defined in $D\times T$ we have
\begin{equation*}
\begin{split}
  &\mathcal{L}\varphi({\bf x},t)\psi({\bf x},t)-\varphi({\bf x},t) \mathcal{L}\psi({\bf x},t)\\
&=\frac{1}{c^2}\partial_t(\partial_t\varphi({\bf x},t) \psi({\bf x},t)-\varphi({\bf x},t) \partial_t \psi({\bf x},t))\\
&-\nabla \cdot(\nabla\varphi({\bf x},t) \psi({\bf x},t)-\varphi({\bf x},t) \nabla \psi({\bf x},t)),
\end{split}
\end{equation*}
so
\begin{align}\label{fundamental3d}
  &\int_{D\times T}\{\mathcal{L}\varphi({\bf x},t)\psi({\bf x},t)-\varphi({\bf x},t) \mathcal{L}\psi({\bf x},t)\}\,\mathrm{d}V\,\mathrm{d}t\nonumber\\
&=\frac{1}{c^2}\int_{D_2}\{\partial_t\varphi({\bf x},t) \psi({\bf x},t)-\varphi({\bf x},t) \partial_t \psi({\bf x},t)\}\,\mathrm{d}V\nonumber\\
&-\frac{1}{c^2}\int_{D_1}\{\partial_t\varphi({\bf x},t) \psi({\bf x},t)-\varphi({\bf x},t) \partial_t \psi({\bf x},t)\}\,\mathrm{d}V\nonumber\\
&-\int_{S\times T}\{\partial_{\bf n}\varphi({\bf x},t) \psi({\bf x},t)-\varphi({\bf x},t) \partial_{\bf n}\psi({\bf x},t)\}\,\mathrm{d}S\,\mathrm{d}t .
\end{align}
This is the fundamental integral identity for the wave equation in 3D.
Next we will need the advanced Green's function for $\mathcal{L}$ which is given by
\begin{equation*}
G({\bf x},t, {\bf x'},t')=Q({\bf x-\bf x'},t-t'),
\end{equation*}
where $Q({\bf y},s)$ satisfies
\begin{equation}\label{eq:3dg}
\frac{1}{c^2}\partial_{ss}Q({\bf y},s)-\nabla^2Q({\bf y},s)=\delta({\bf y})\delta(s),\qquad s<0,
\end{equation}
and 
$$Q({\bf y},s)=0,\qquad s>0.$$
Because of translational invariance, we can take the Fourier transform of (\ref{eq:3dg}) and get
\begin{equation*}
[-(\frac{\omega}{c})^2+\zeta^2]\hat Q({\bf \xi},\omega)=1,\qquad  \zeta=|{\bf \xi} |,
\end{equation*}
\begin{equation}\label{Qequation}
\hat Q({\bf \xi},\omega)=\frac{1}{D({\bf \xi},\omega)},
\end{equation}
where
\begin{equation*}
D({\bf \xi},\omega)=-(\frac{\omega}{c})^2+\zeta^2.
\end{equation*}
Applying the inverse Fourier transform on (\ref{Qequation}) gives
\begin{equation*}
Q({\bf y},s)=\frac{1}{16\pi^4}\int q({\bf \xi}, s) e^{i {\bf \xi} \cdot {\bf y}}\,\mathrm{d}{\bf \xi},
\end{equation*}
where
\begin{equation*}
q({\bf \xi},s)=\int_{C_\varepsilon} \frac{e^{-izs}}{D({\bf \xi},z)}\,\mathrm{d}z,
\end{equation*}
and $C_\varepsilon$ is a contour slightly below the real axis. See figure \ref{contour}. 
Observe that the integrand has simple poles at $z=\pm c\zeta$, so if $s>0$, we must close the contour in the lower half-plane. By Cauchy, $$q({\bf \xi},s)=0.$$ If $s<0,$ we must now close the contour in the upper half plane and Cauchy's theorem gives
\begin{equation*}
\begin{split}
q({\bf \xi},s)=\frac{c\pi i}{\zeta}[e^{isc\zeta}-e^{-isc\zeta}],
\end{split}
\end{equation*}
and
\begin{equation*}
\begin{split}
\frac{1}{16\pi^4}\int \,\mathrm{d} {\bf \xi}\ \frac{c\pi i}{\zeta}e^{isc\zeta}e^{i{\bf \xi} \cdot {\bf y}}=\frac{c}{8\pi^2 r}[\int_0^{\infty} \,\mathrm{d} \zeta\ e^{i(r+sc)\zeta}-\int_0^{\infty} \,\mathrm{d} \zeta\ e^{-i(r-sc)\zeta}],
\end{split}
\end{equation*}
where $ r=|\bf y|.$
Similarly,
\begin{equation*}
-\frac{1}{16\pi^4}\int \,\mathrm{d} {\bf \xi}\frac{c\pi i}{\zeta}e^{-isc\zeta}e^{i{\bf \xi} \cdot {\bf y}}=\frac{-c}{8\pi^2 r}[\int_0^{\infty} \,\mathrm{d} \zeta\ e^{i(r-sc)\zeta}-\int_0^{\infty} \,\mathrm{d} \zeta\ e^{-i(r+sc)\zeta}].
\end{equation*}
In the end, 
\begin{equation*}
\begin{split}
Q({\bf y},s)=\frac{c}{4\pi r}[\delta(r+sc)-\delta(r-sc)].
\end{split}
\end{equation*}
Since if $s<0,$ $\delta(r-sc)=0,$
finally
\begin{equation*}
Q({\bf y},s)=\frac{1}{4\pi r}\delta(s+\frac{r}{c}).
\end{equation*}
Thus our advanced Green's function is
\begin{equation}\label{3dgreen}
G({\bf x},t,{\bf x'},t') = \left\{
  \begin{array}{l l}
    0 & \quad \text{if $t>t'$ }\\
    \frac{1}{4\pi |{\bf x-\bf x'}|}\delta(t-t'+\frac{|{\bf x-\bf x'}|}{c}), & \quad \text{if $t<t'$}.
  \end{array} \right.
\end{equation}
We now apply (\ref{3dgreen}) to the fundamental integral identity (\ref{fundamental3d}). Let $\varphi({\bf x},t)$ be a solution to the wave equation (\ref{3dwaveequation})
$$
\mathcal{L}\varphi({\bf x},t)=\rho({\bf x},t),
$$
then for any $({\bf x},t) \in D\times T$, we have
\begin{equation*}
\begin{split}
&\int_{D\times T} \,\mathrm{d}V'\,\mathrm{d}t'\{\mathcal{L}\varphi({\bf x'},t')G({\bf x'},t',{\bf x},t)-\varphi({\bf x'},t')\mathcal{L}G({\bf x'},t',{\bf x},t)\}\\
&=\frac{1}{c^2}\int_{D_2} \,\mathrm{d}V'\{\partial_{t'}\varphi({\bf x'},t')G({\bf x'},t',{\bf x},t)-\varphi({\bf x'},t')\partial_{t'}G({\bf x'},t',{\bf x},t)\}|_{t'=t_2}\\
&-\frac{1}{c^2}\int_{D_1} \,\mathrm{d}V'\{\partial_{t'}\varphi({\bf x'},t')G({\bf x'},t',{\bf x},t)-\varphi({\bf x'},t')\partial_{t'}G({\bf x'},t',{\bf x},t)\}|_{t'=t_1}\\
&-\int_{S\times T} \,\mathrm{d}S'\,\mathrm{d}t'\{\partial_{{\bf n}'}\varphi({\bf x'},t')G({\bf x'},t',{\bf x},t)-\varphi({\bf x'},t')\partial_{{\bf n}'}G({\bf x'},t',{\bf x},t)\}.
\end{split}
\end{equation*}
Since $t_2>t,$ for the advance Green's function $$G({\bf x'},t',{\bf x},t)|_{t'=t_2}=0.$$  And since the field is entirely driven by the source $$\varphi({\bf x},t_1)=0,$$ we get
\begin{equation*}
\begin{split}
\varphi({\bf x},t)&=\int_{D\times T}\,\mathrm{d}V'\,\mathrm{d}t'\rho({\bf x'},t')G({\bf x'},t',{\bf x},t)\\
&+\int_{S\times T} \,\mathrm{d}S'\,\mathrm{d}t'\{\partial_{{\bf n}'}\varphi({\bf x'},t')G({\bf x'},t',{\bf x},t)-\varphi({\bf x'},t')\partial_{{\bf n}'}G({\bf x'},t',{\bf x},t)\}.
\end{split}
\end{equation*}
and
\begin{equation*}
G({\bf x'},t',{\bf x},t)=h({\bf x'},{\bf x})\theta(t-t') \delta(|{\bf x'}-{\bf x}|+c(t'-t)),
\end{equation*}
where $$
h({\bf x'},{\bf x})=\frac{c}{4\pi |{\bf x'}-{\bf x}|}.
$$
Since
\begin{equation*}
\begin{split}
\partial_{{\bf n}'}G({\bf x'},t',{\bf x},t)&=\theta(t-t')\partial_{{\bf n}'}h({\bf x'},{\bf x})\delta(|{\bf x'}-{\bf x}|+c(t'-t))\\
&+\theta(t-t')h({\bf x'},{\bf x})\partial_{{\bf n}'}\delta(|{\bf x'}-{\bf x}|+c(t'-t)),
\end{split}
\end{equation*}
and
\begin{equation*}
\begin{split}
\partial_{{\bf n}'}\delta(|{\bf x'}-{\bf x}|+c(t'-t))=\frac{1}{c}\partial_{{\bf n}'}|{\bf x'}-{\bf x}|\partial_{t'}\delta(|{\bf x'}-{\bf x}|+c(t'-t)),
\end{split}
\end{equation*}
 we thus have
\begin{equation*}
\begin{split}
\int_{D\times T}\,\mathrm{d}V'\,\mathrm{d}t'\rho({\bf x'},t')G({\bf x'},t',{\bf x},t)=\int_D\,\mathrm{d}V'h({\bf x'},{\bf x})\rho({\bf x'},T),
\end{split}
\end{equation*}
where
\begin{equation*}
 T=T(t,{\bf x'},{\bf x})=t-\frac{1}{c}|{\bf x'}-{\bf x}|.
\end{equation*}
From
\begin{equation*}
\int_{S\times T}\,\mathrm{d}S'\,\mathrm{d}t'\partial_{{\bf n}'}\varphi({\bf x'},t')G({\bf x'},t',{\bf x},t)=\int_S\, \mathrm{d} S' h({\bf x'},{\bf x})(\partial_{{\bf n}'}\varphi)({\bf x'},T),
\end{equation*}
and
\begin{equation*}
\begin{split}
&-\int_{S\times T}\,\mathrm{d}S'\,\mathrm{d}t'\varphi({\bf x'},t')\partial_{{\bf n}'}G({\bf x'},t',{\bf x},t)\\
&=-\int_S\,\mathrm{d}S'\partial_{{\bf n}'}h({\bf x'},{\bf x})\varphi({\bf x'},T)+\int_{S}\,\mathrm{d}S'\frac{1}{c}h({\bf x'},{\bf x})\partial_{{\bf n}'}|{\bf x'}-{\bf x}|(\partial_{t'}\varphi)({\bf x'},T).
\end{split}
\end{equation*}
finally we get
\begin{equation}\label{eq:mfundamental}
\begin{split}
\varphi({\bf x},t)&=\int_D\,\mathrm{d}V' h({\bf x'},{\bf x})\rho({\bf x'},T)+\int_S\,\mathrm{d}S'\{h({\bf x'},{\bf x})(\partial_{{\bf n}'}\varphi)({\bf x'},T)\\
&-\partial_{{\bf n}'}h({\bf x'},{\bf x})\varphi({\bf x'},T)+\frac{h({\bf x'},{\bf x})}{c}\partial_{{\bf n}'}|{\bf x'}-{\bf x}|(\partial_{t'}\varphi)({\bf x'},T)\}.
\end{split}
\end{equation}
This is the integral identity for an operator defining a 3D
wave equation and it holds for any solution of the scalar 3D wave equation.

\section{The integral identity of the electric wave equation}\label{Aintegral}
Here we do some calculations to derive  (\ref{esolution}) from (\ref{eq:original e solution}).  For the writing in simplicity, we write ${\bf E}, {\bf J}, \rho, c, \mu, \epsilon$ instead of ${\bf E}_j, {\bf J}_j, \rho_j, c_j, \mu_j, \epsilon_j,$ $j=0,1$ respectively here. Observe first that
\begin{equation*}
\begin{split}
\partial_{{\bf n}'}({\bf E}({\bf x'},T))&=({\bf n'}\cdot \nabla')({\bf E}({\bf x'},T))=(({\bf n'}\cdot \nabla'){\bf E})({\bf x'},T)\\
&+(\partial_{t'}E)({\bf x'},T)(-\frac{1}{c}({\bf n'}\cdot \nabla')|{\bf x'}-{\bf x}|),
\end{split}
\end{equation*}
so,
\begin{equation}\label{eq:original e solution 2}
\begin{split}
{\bf E}({\bf x},t)&=-\int_D\,\mathrm{d}V'h({\bf x'},{\bf x})\{\mu \partial_{t'}{\bf J}+\frac{1}{\varepsilon}\nabla'\rho\}({\bf x'},T)\\
&+\int_S\,\mathrm{d}S'\{h({\bf x'},{\bf x})\partial_{{\bf n}'}({\bf E}({\bf x'},T))+\frac{1}{c}h({\bf x'},{\bf x})(\partial_{t'}{\bf E})({\bf x'},T)\partial_{{\bf n}'}|{\bf x'}-{\bf x}|\\&-\partial_{{\bf n}'}h({\bf x'},{\bf x}){\bf E}({\bf x'},T)+\frac{1}{c}h({\bf x'},{\bf x})(\partial_{t'}{\bf E})({\bf x'},T)\partial_{{\bf n}'}|{\bf x'}-{\bf x}|\}\\
&+\int_S\,\mathrm{d}S'\{h({\bf x'},{\bf x})\partial_{{\bf n}'}({\bf E}({\bf x'},T))-\partial_{{\bf n}'}h({\bf x'},{\bf x}){\bf E}({\bf x'},T)\\
&+\frac{2}{c}h({\bf x'},{\bf x})(\partial_{t'}{\bf E})({\bf x'},T)\partial_{{\bf n}'}|{\bf x'}-{\bf x}|\}.
\end{split}
\end{equation}
We are going to rework the first term in the integral (\ref{eq:original e solution 2}). Observe that for a vector field ${\bf a}$ and a scalar $f$ we have,
\begin{equation*}
\begin{split}
({\bf n} \cdot \nabla)(f {\bf a})&=({\bf n} \cdot \nabla f){\bf a}+f({\bf n} \cdot \nabla){\bf a},\\
\nabla \cdot (f{\bf a})&=\nabla f \cdot {\bf a}+f\nabla \cdot {\bf a},\\
\nabla \times (f{\bf a})&=\nabla f\times {\bf a}+f\nabla \times {\bf a},
\end{split}
\end{equation*}
so,
\begin{equation*}
\begin{split}
{\bf n} \times (\nabla \times (f{\bf a}))=\nabla f({\bf n} \cdot {\bf a})-{\bf a}({\bf n}\cdot \nabla f )+f{\bf n}\times(\nabla \times {\bf a}).
\end{split}
\end{equation*}
Further,
\begin{equation*}
\begin{split}
&({\bf n} \cdot \nabla)(f {\bf a})+{\bf n} \times (\nabla \times (f{\bf a}))-n\nabla\cdot (f{\bf a})\\
&=f({\bf n} \cdot \nabla){\bf a}+({\bf n} \times\nabla f)\times{\bf a}+f{\bf n}\times(\nabla \times  {\bf a})-f{\bf n}\nabla\cdot{\bf a}.
\end{split}
\end{equation*}
so if we let $f=h({\bf x'},{\bf x})$ and ${\bf a}={\bf E}({\bf x'},T), $ we thus have
\begin{equation}\label{equation1}
\begin{split}
&\int_S\,\mathrm{d}S' h({\bf x'},{\bf x})\partial_{{\bf n}'}({\bf E}({\bf x'},T))\\
&=\int_S\,\mathrm{d}S'\{
-({\bf n'}\times\nabla'h({\bf x'},{\bf x}))\times{\bf E}({\bf x'},T)\\
&-h({\bf x'},{\bf x}){\bf n'}\times(\nabla'\times({\bf E}({\bf x'},T)))\\
&+h({\bf x'},{\bf x}){\bf n'}\nabla'\cdot({\bf E}({\bf x'},T))\}.
\end{split}
\end{equation}
Inserting (\ref{equation1}) into (\ref{eq:original e solution 2}) leads to
\begin{equation*}
\begin{split}
{\bf E}({\bf x},t)&=-\int_D\,\mathrm{d}V' h({\bf x'},{\bf x})\{\mu \partial_{t'}{\bf} J+\frac{1}{\varepsilon}\nabla'\rho\}({\bf x},T)\\
&+\int_S\,\mathrm{d}S'\{
h({\bf x'},{\bf x}){\bf n'}\nabla'\cdot({\bf E}({\bf x'},T))-({\bf n'}\times\nabla'h({\bf x'},{\bf x}))\times {\bf E}({\bf x'},T)\\
&-h({\bf x'},{\bf x}){\bf n'}\times(\nabla'\times({\bf E}({\bf x'},T)))-\partial_{{\bf n}'}h({\bf x'},{\bf x}){\bf E}({\bf x'},T)\\
&+\frac{2}{c}h({\bf x'},{\bf x})(\partial_{t'}{\bf E})({\bf x},T)\partial_{{\bf n}'}|{\bf x'}-{\bf x}|\}.
\end{split}
\end{equation*}
Since
\begin{equation*}
\begin{split}
\nabla'({\bf E}({\bf x'},T))&=(\nabla'\cdot {\bf E})({\bf x'},T)-\frac{1}{c}(\partial_{t'} {\bf E})({\bf x'},T)\cdot\nabla'|{\bf x'}-{\bf x}|,\\
\nabla'\times({\bf E}({\bf x'},T))&=(\nabla'\times {\bf E})({\bf x'},T)+\frac{1}{c}(\partial_{t'} {\bf E})({\bf x'},T)\times\nabla'|{\bf x'}-{\bf x}|,
\end{split}
\end{equation*}
in the end, (\ref{eq:original e solution 2}) can be written in the following form
\begin{equation*}
\begin{split}
{\bf E}({\bf x},t)&=-\int_D\,\mathrm{d}V' h({\bf x'},{\bf x})\{\mu \partial_{t'}{\bf J}+\frac{1}{\varepsilon}\nabla'\rho\}({\bf x'},T)\\
&+\int_S\,\mathrm{d}S'\{h({\bf x'},{\bf x}){\bf n'}(\nabla'\cdot {\bf E})({\bf x'},T)\\
&-\frac{1}{c}h({\bf x'},{\bf x})\bf n'(\partial_{t'} {\bf E})({\bf x'},T)\cdot\nabla'|{\bf x'}-{\bf x}|)\\
&-({\bf n'}\times\nabla'h({\bf x'},{\bf x}))\times {\bf E}({\bf x'},T)-h({\bf x'},{\bf x}){\bf n'}\times (\nabla'\times{\bf E})({\bf x'},T)\\
&-\frac{1}{c}h({\bf x'},{\bf x}){\bf n'}\times((\partial_{t'}{\bf E})({\bf x'},T)\times\nabla'|{\bf x'}-{\bf x}|)\\
&-\partial_{{\bf n}'}h({\bf x'},{\bf x}){\bf E}({\bf x'},T)+\frac{2}{c}h({\bf x'},{\bf x})(\partial_{t'}{\bf E})({\bf x'},T)\partial_{{\bf n}'}|{\bf x'}-{\bf x}|\}.
\end{split}
\end{equation*}
Notice that
\begin{equation*}
\begin{split}
&({\bf n}\times{\bf a})\times \nabla f-({\bf n}\times\nabla f)\times{\bf a}=({\bf n}\cdot\nabla f){\bf a}-({\bf n}\cdot {\bf a})\nabla f,
\end{split}
\end{equation*}
and
\begin{equation*}
-({\bf n}\times\nabla f)\times{\bf a}-({\bf n}\cdot\nabla f){\bf a}=-({\bf n}\times{\bf a})\times \nabla f-({\bf n}\cdot {\bf a})\nabla f,
\end{equation*}
and performing them on $h$ and ${\bf E}$ gives
\begin{equation*}
\begin{split}
&-({\bf n'}\times\nabla' h({\bf x'},{\bf x}))\times{\bf E}({\bf x'},T)-\partial_{{\bf n}'}h({\bf x'},{\bf x}){\bf E}({\bf x'},T)\\
&=-({\bf n'}\times{\bf E}({\bf x'},T))\times \nabla' h({\bf x'},{\bf x})-({\bf n'}\cdot {\bf E}({\bf x'},T))\nabla' h({\bf x'},{\bf x}).
\end{split}
\end{equation*}
In addition
\begin{equation*}
\begin{split}
&({\bf n}\times{\bf a})\times\nabla f+{\bf n}\times({\bf a} \times\nabla f)\\
&=2({\bf n}\cdot \nabla f){\bf a}-({\bf a}\cdot\nabla f){\bf n}-({\bf n}\cdot{\bf a})\nabla f,
\end{split}
\end{equation*}
and
\begin{equation*}
\begin{split}
&-({\bf a}\cdot\nabla f){\bf n}-{\bf n}\times({\bf a}\times\nabla f)+2({\bf n}\cdot \nabla f){\bf a}\\
&=({\bf n}\times{\bf a})\times\nabla f+({\bf n}\cdot{\bf a})\nabla f,
\end{split}
\end{equation*}
give
\begin{equation*}
\begin{split}
&-{\bf n'}((\partial_{t'} {\bf E})({\bf x'},T)\cdot\nabla'|{\bf x'}-{\bf x}|)-{\bf n'}\times ((\partial_{t'}{\bf E})({\bf x},T)\times\nabla'|{\bf x'}-{\bf x}|)\\
&+2\partial_{{\bf n}'}|{\bf x'}-{\bf x}|(\partial_{t'}{\bf E})({\bf x'},T)\\
&=({\bf n'}\times(\partial_{t'}{\bf E})({\bf x'},T))\times\nabla' |{\bf x'}-{\bf x}|+({\bf n'}\cdot(\partial_{t'}{\bf E})({\bf x'},T))\nabla' |{\bf x'}-{\bf x}|.
\end{split}
\end{equation*}
Thus
\begin{equation*}
\begin{split}
{\bf E}({\bf x},t)&=-\int_D\,\mathrm{d}V'h({\bf x'},{\bf x})\{\mu \partial_{t'}{\bf J}+\frac{1}{\varepsilon}\nabla'\rho\}({\bf x'},T)\\
&+\int_S\,\mathrm{d}S'\{h({\bf x'},{\bf x}){\bf n'}(\nabla'\cdot {\bf E})({\bf x'},T)-({\bf n'}\times E({\bf x'},T))\times\nabla'h({\bf x'},{\bf x})\\
&-({\bf n'}\cdot E({\bf x'},T))\nabla'h({\bf x'},{\bf x})+\frac{1}{c}h({\bf x'},{\bf x})({\bf n'}\times (\partial_{t'}{\bf E})({\bf x'},T))\times \nabla'|{\bf x'}-{\bf x}|\\
&+\frac{1}{c}h({\bf x'},{\bf x})({\bf n'}\cdot (\partial_{t'}{\bf E})({\bf x'},T))\nabla'|{\bf x'}-{\bf x}|-h({\bf x'},{\bf x}){\bf n'}\times(\nabla'\times {\bf E})({\bf x'},T)\}.
\end{split}
\end{equation*}
Using the special form of the divergence theorem, we have
\begin{equation*}
\begin{split}
&\int_S\,\mathrm{d}S' h({\bf x'},{\bf x}){\bf n'}(\nabla'\cdot {\bf E})({\bf x'},T)\\
&=\int_D\,\mathrm{d}V' h({\bf x'},{\bf x})\frac{1}{\varepsilon}(\nabla'\rho)({\bf x'},T)-\int_D,\mathrm{d}V' h({\bf x'},{\bf x})\frac{1}{\varepsilon}\nabla\rho({\bf x'},T)\\
&+\int_D\,\mathrm{d}V'\frac{1}{\varepsilon}\rho({\bf x'},T)\nabla'h({\bf x'},{\bf x}),
\end{split}
\end{equation*}
where
\begin{equation*}
\begin{split}
\nabla'h({\bf x'},{\bf x})=\frac{1}{4\pi}\nabla'\frac{1}{|{\bf x'}-{\bf x}|}=-\nabla h({\bf x'},{\bf x}).
\end{split}
\end{equation*}
Together with
\begin{equation*}
(\partial_{t'}{\bf J})({\bf x'},T)=\partial_t({\bf J}({\bf x'},T)),
\end{equation*}
we finally get
\begin{equation*}
\begin{split}
&{\bf E}({\bf x},t)=-\partial_t\frac{\mu}{4\pi}\int_D\,\mathrm{d}V'\frac{{\bf J}({\bf x'},T)}{|{\bf x'}-{\bf x}|}-\nabla\frac{1}{4\pi\varepsilon}\int_D\,\mathrm{d}V'\frac{\rho({\bf x'},T)}{|{\bf x'}-{\bf x}|}\\
&+\partial_t[\frac{1}{4\pi}\int_S\,\mathrm{d}S'\{\frac{1}{c|{\bf x'}-{\bf x}|} ({\bf n'}\times {\bf E}({\bf x'},T))\times \nabla'|{\bf x'}-{\bf x}|\\
&+\frac{1}{c|{\bf x'}-{\bf x}|} ({\bf n'}\cdot {\bf E}({\bf x'},T))\nabla'|{\bf x'}-{\bf x}|+\frac{1}{|{\bf x'}-{\bf x}|} {\bf n'}\times {\bf B}({\bf x'},T)  \}]\\
&-\frac{1}{4\pi}\int_S\,\mathrm{d}S'\{ ({\bf n'}\times {\bf E}({\bf x'},T))\times \nabla'\frac{1}{|{\bf x'}-{\bf x}|}\\
&+({\bf n'}\cdot {\bf E}({\bf x'},T))\nabla'\frac{1}{|{\bf x'}-{\bf x}|}\}.
\end{split}
\end{equation*}
This is the integral identity of the electric  wave equation (\ref{ewave}).

\end{appendices}

\bibliographystyle{unsrt}
\bibliography{paper2}

\begin{thebibliography}{10}

\bibitem{PML1}
J.~Berenger.
\newblock A perfectly matched layer for the absorption of electromagnetic
  waves.
\newblock {\em Journal of computational physics}, 114(2):185--200, 1994.

\bibitem{PML2}
S.~Gedney.
\newblock An anisotropic perfectly matched layer absorbing media for the
  truncation of fdtd lattices.
\newblock {\em IEEE Transactions on Antennas and Propagation},
  44(12):1630--1639, 1996.

\bibitem{FEM1}
S.~Ahmed.
\newblock Finite-element method for waveguide problems.
\newblock {\em Electronics Letters}, 4(18):387--389, 1968.

\bibitem{FEM2}
J.~P. Webb.
\newblock Application of the finite-element method to electromagnetic and
  electric topics.
\newblock {\em Rep. Prog. Phys.}, 58:1673--1712, 1995.

\bibitem{FEM3}
J.~Jin.
\newblock {\em The Finite Element Method in Electromagnetics, 2nd Edition}.
\newblock Wiley-IEEE Press, 2002.

\bibitem{Cheng}
Alexander H.~H. Cheng and Daisy~T. Cheng.
\newblock Heritage and early history of the boundary element method.
\newblock {\em {Engineering Vnalysis with Boundary Elements}}, 29:268--302,
  2005.

\bibitem{Jones64}
D.~S. Jones.
\newblock {\em The Theory of Electromagnetism}, volume~47 of {\em International
  Series of Monographs on Pure and Applied mathematics}.
\newblock Pergamon Press, 1964.

\bibitem{Weile}
Greeshma Pisharody and Daniel~S. Weile.
\newblock Electromagnetic scattering from homogeneous dielectric bodies using
  time-domain integral equations.
\newblock {\em IEEE Transactions on Antennas and Propagation}, 54:687--697,
  2006.

\bibitem{Weile2}
D.~S. Weile, I.~Uluer, J.~Li, and D.A. Hopkins.
\newblock Integration rules and experimental evidences for the stability of
  time domain integral equations.
\newblock {\em International Applied Computational Electromagnetics Society
  Symposium(ACES)}, 2017.

\bibitem{Shanker}
A.~J. Pray, N.~V. Nair, and B.~Shanker.
\newblock Stability properties of the time domain electric field integral
  equationusing a separable approximation for the convolution with the retarded
  potential.
\newblock {\em IEEE Transactions on Antennas and Propagation}, 60:3772--3781,
  2012.

\bibitem{Walker}
M.~J. Bluck and S.~P. Walker.
\newblock Time-domain bie analysis of large three-dimensional electromagnetic
  scattering problems.
\newblock {\em IEEE Transactions on Antennas and Propagation}, 45:894--901,
  1997.

\bibitem{Walker2}
S.~Dodson, S.~P. Walker, and M.~J. Bluck.
\newblock Implicitness and stability of time domain integral equation
  scattering analysis.
\newblock {\em Appl. Comput. Electromagn. Soc. J.}, 13:291--301, 1998.

\bibitem{Zhao}
Ying Zhao, Dazhi Ding, and Rushan Chen.
\newblock A discontinuous galerkin time-domain integral equation method for
  electromagnetic scattering from pec objects.
\newblock {\em IEEE Transactions on Antennas and Propagation},
  64(6):2410--2415, 2016.

\bibitem{Huang}
Li~Huang, Yi-Bei Hou, Hao-Xuan Zhang, Liang Zhou, and Wen-Yan Yin.
\newblock A discontinuous galerkin time-domain integral equation method for
  electromagnetic scattering from pec objects.
\newblock {\em IEEE Transactions on Electromagnetic Compatibility}, 2018.

\bibitem{wolf}
D.~N. Pattanayak and E.~Wolf.
\newblock General form and new interpretation of the ewald-oseen extincition
  theorem.
\newblock {\em Optics communications}, 6(3):217--220, 1972.

\bibitem{Aihua}
Aihua Lin, Anastasiia Kuzmina, and Per~Kristen Jakobsen.
\newblock A boundary integral approach to linear and nonlinear transient wave
  scattering.
\newblock {\em Submitted}, 2018.

\bibitem{Jakobsen09}
Per Jakobsen.
\newblock Calculating optical forces using the boundary integral method.
\newblock {\em Physics Scripta}, 80:35401--9, 2009.

\bibitem{Aihua3}
Aihua Lin and Per~Kristen Jakobsen.
\newblock Numerical techniques of solving a 3d nonlinear maxwell’s equations
  by the eos formulations: Parallelization, stabilities and singularities.
\newblock {\em To appear}, 2018.

\end{thebibliography}

\end{document}